\documentclass{article}
\begin{document}
\input amssym.def

\def \C{\Bbb C}
\def \d{{\partial}}
\def \g{\gamma}
\def \G{\Gamma}
\def \H{\Bbb H}
\def \l{\lambda}
\def \L{\Lambda}
\def \la{\langle}
\def \O{\Omega}
\def \P{\Bbb P}
\def \ra{\rangle}
\def \N{\Bbb Nll}
\def \R{\Bbb R}
\def \sm{\setminus}
\def \summa{\sum_{\g\in\G_0\sm\G}}
\def \w{\wedge}
\def \ww{\la w,w\ra}
\def \xy{\la X,Y\ra}
\def \Z{\Bbb Z}
\def \zz{\la z,z\ra}
\def \zx{\la z,X\ra}
\def \zy{\la z,Y\ra}
\def \elt{\pmatrix{
a_{11} & \dots & a_{1n} & b_1 \cr
\dots & \dots & \dots & \dots \cr
a_{n1} & \dots & a_{nn} & b_n \cr
c_1 &  \dots & c_n & d
}
}

\newtheorem{Th}{Theorem}[section]
\newtheorem{assum}[Th]{Assumption}
\newtheorem{cor}[Th]{Corollary}
\newtheorem{definition}[Th]{Definition}
\newtheorem{lem}[Th]{Lemma}
\newtheorem{prop}[Th]{Proposition}
\newtheorem{remark}[Th]{Remark}

\def \proof{{\noindent{\it Proof.\ \ }}}

\author{Tatyana Foth}

\title{Bohr-Sommerfeld tori and relative Poincar\'e series on a complex hyperbolic
space }

\maketitle

\abstract{Automorphic forms on a bounded symmetric domain $D=G/K$ can be
viewed as holomorphic sections
of $L^{\otimes k}$, where $L$ is a quantizing line bundle on a compact quotient of $D$
and $k$ is a positive integer.

Let $\G$ be a cocompact discrete subgroup of $SU(n,1)$ which acts freely on $SU(n,1)/U(n)$. 
We suggest a construction of 
relative Poincar\'e series associated to loxodromic elements of $\G$.
In complex dimension 2 we describe the Bohr-Sommerfeld tori in 
$\G\backslash SU(n,1)/U(n)$ associated to hyperbolic elements of $\G$  
and prove that the relative Poincar\'e series associated to hyperbolic 
elements of $\G$ are not identically zero for large values of $k$. }

\section{Introduction}

\subsection{General definitions}

We shall start with a brief review of the general concept of an 
automorphic form. Let $G$ be a connected non-compact real 
semi-simple Lie group,
$K$ be a maximal compact subgroup of $G$, $\G$ be a discrete 
subgroup of $G$ such that $\G\backslash G$ has a finite volume. 
Let $V$ be a finite-dimensional vector space, $\rho:K\rightarrow GL(V)$
be a representation of K. A smooth $Z({\frak g})$-finite function 
$f:G\rightarrow V$ is called an {\bf automorphic form on $\bf G$ for $\bf \G$} 
if 
\begin{equation}
f(\g gk)=f(g)\rho(k)
\label{aut-law}
\end{equation}
for any $\g\in \G$, $g\in G$, $k\in K$, and there are a positive constant
$C$ and a non-negative integer $m$ such that  
\begin{equation}
|f(g)|\le C||g||^m
\label{growth}
\end{equation}   
for any $g\in G$, here $|.|$ is the norm corresponding to 
a $\rho(K)$-invariant Hilbert structure on $V$, $||g||=tr (g^*g)$
taken in the adjoint representation of $G$. An automorphic form $f$ 
is called a {\it cusp} form if $f\in L^\infty(\G\backslash G)$.

The automorphy law (\ref{aut-law}) means geometrically
that $f$ defines a $\G$-invariant section of the vector bundle
$G\times _K V\rightarrow G/K$ associated to the principal bundle
$G\rightarrow G/K$, here $G\times _K V=G\times V/\sim$, and
the equivalence relation is given by the representation $\rho$:
$(g,v)\sim (gk,v\rho(k))$. 

The growth condition (\ref{growth}) 
is automatically satisfied with $m=0$ in the case when $\G\backslash G$ 
is compact and in this case any automorphic form is a cusp form.

Recall also that a function $f:G\rightarrow V$ is said to be 
{\bf ${\bf Z({\frak g})}$-finite} if it is annihilated by an ideal $I$ 
of $Z({\frak g})$ of a finite codimension, here $Z({\frak g})$ is the center
of the universal enveloping algebra $U({\frak g})$. 

$U({\frak g})$ can be identified with
the algebra $D(G)$ of all left-invariant differential operators
on $G$: to $Y\in \frak g$ is associated a differential operator
$Yf(g)=\frac{d}{dt}f(ge^{tY})|_{t=0}$, this establishes a linear map
${\frak g}\rightarrow D(G)$ which extends to an isomorphism $U({\frak g})
\rightarrow D(G)$. $Z({\frak g})$ can be viewed as the subalgebra 
of all bi-invariant differential operators, it is isomorphic
to a polynomial ring in $l$ letters where $l$ is the rank of $G$.
A useful example to have in mind is $G=SL(2,R)$ and $codim \ I=1$, then 
we have: $l=1$, $Z({\frak g})$ is generated by the Casimir operator $\cal C$,
and saying that a function $f$ is $Z({\frak g})$-finite is equivalent
to stating that $f$ is an eigenfunction of $\cal  C$.

A well-known construction of an automorphic form on $G$ is {\bf Poincar\'e
series}
$$
\sum_{\g\in\G}q(\g g),
$$
where the function $q:G\rightarrow V$ is $Z({\frak g})$-finite and 
{\it $K$-finite on the right} 
(i.e. the set of its right translates under elements of $K$ is a 
finite-dimensional vector space), and $q\in L^1(G)$.
One can also consider {\bf relative Poincar\'e series}
$$
\sum_{\g\in\G_0\backslash\G}q(\g g),
$$
where $q:G\rightarrow V$ is $Z({\frak g})$-finite, 
$K$-finite on the right, $\G_0$-invariant, 
and $q\in L^1(\G_0\backslash G)$. 

Let us explain now how to construct an automorphic form on $G/K$. 
An {\bf automorphy factor} is a map $\mu:\G\times G/K\rightarrow GL(V)$
such that $\mu(g_1g_2,x)=\mu(g_1,g_2x)\mu(g_2,x)$. It allows to
define an {\bf automorphic form on $\bf G/K$} as a function 
$f:G/K\rightarrow V$ such that 
$$
f(\g x)\mu(\g,x)=f(x)
$$ 
for any
$\g\in\G$, $x\in G/K$. Notice that then the function $F(g)=f(g(0))\mu(g,0)$,
where $g\in G$, $x=g(0)\in G/K$, satisfies (\ref{aut-law}) 
with $\rho(k)=\mu(k,0)$. Here $0$ is the fixed point of $K$ in $G/K$.
If $f$ is holomorphic then $F$ is $Z({\frak g})$-finite.

In particular, for a smooth function $q\in L^1(G/K)$
the Poincar\'e series on $G/K$ is
\begin{equation}
\sum_{\g\in\G}q(\g x)\mu(\g,x).
\label{p-s}
\end{equation}

Similarly for a smooth $\G_0$-invariant function 
$q\in L^1(\G_0\backslash G/K)$ the relative Poincar\'e series is 
\begin{equation}
\sum_{\g\in\G_0\backslash \G}q(\g x)\mu(\g,x).
\label{rp-s}
\end{equation}

\subsection{Automorphic forms on bounded symmetric domains
and quantization}

Consider a classical system $(M,\omega)$, where $M$ is a manifold,
and $\omega$ is a symplectic form on $M$. 
The main problem of quantization is to associate a quantum system
$({\cal H}, {\cal O})$ to $(M,\omega)$, where $\cal  H$ is a Hilbert space
and $\cal  O$ consists of symmetric operators on $\cal  H$. 
 

The map $f\mapsto \hat f$,
where $f\in C^\infty(M)$
and $\hat f\in \cal O$, should satisfy the following requirements: 
1) it is $\Bbb R$-linear,
2) if $f=const$ then $\hat f$ is the corresponding multiplication
operator,
3) if $\{ f_1,f_2 \}=f
_3$ then $\hat f_1\hat f_2-\hat f_2\hat f_1=-i\hbar
\hat f_3$.  


How do automorphic forms appear in the context of quantization ?

Suppose that $M$ is a compact K\"ahler manifold which is a quotient
of a bounded symmetric domain $D=G/K$ by the action of a discrete
subgroup $\G$, i.e. $M=\G\backslash D$. Then the quantum phase space
$\cal H$ consists of holomorphic automorphic forms on $D$ for $\G$.
More precisely, let us consider the well-known quantization scheme 
for compact K\"ahler manifolds via Toeplitz operators
(it is related to the standard scheme of geometric quantization
with K\"ahler polarization). Then 
automorphic forms are holomorphic sections of $L^{\otimes k}$, where
the canonical line bundle $L=\Lambda ^nT^*M$ 
is the quantizing line bundle on $M$, here $n=\dim_\C M$ and $k$
is a positive integer which determines the weight of an automorphic
form, and $\hbar=\frac{1}{k}$. 



We also notice that the automorphic forms 
(\ref{p-s}) and (\ref{rp-s}) 
are sums of coherent states associated to holomorphic discrete
series representations of $G$.

Let us describe all this in a bit more details. 
Let $D=G/K$ be a bounded symmetric domain, it is a Hermitian symmetric
space of noncompact type (so the Riemannian metric on $D$ is given
by the real part of the hermitian form, and the symplectic form,
which is a K\"ahler form in this case, is given by the imaginary part
of the hermitian form, all these forms are $G$-invariant, of course). 
The irreducible Hermitian spaces of non-compact type are 

\noindent I) $SU(p,q)/S(U(p)\times U(q))$,

\noindent II) $Sp(p,\R)/U(p)$,
 
\noindent III) $SO^*(2p)/U(p)$,

\noindent IV) $SO_o(p,2)/SO(p)\times SO(2)$ 

\noindent (and let us omit the case of an exceptional Lie algebra). 
So we have a metric 
\begin{equation}
ds^2=g_{ij}dz^id\bar z^j,
\label{metric}
\end{equation}
the corresponding K\"ahler form is $\omega =g_{ij}dz^i\w d\bar z^j=i
\d\bar\d 
\ln K(z,z)$, where $K(z,w)$ is the Bergman kernel 
of the domain $D$. Recall that $K(z,w)=\overline{K(w,z)}$ and
$$
K(\g z,\g w)=[\det J(\g,z)]^{-1}[\det \bar J(\g,w)]^{-1}K(z,w).
$$ 
The Poisson bracket is
$$
\{ f,g\}=ig^{jl}(\frac{\d f}{\d z^j}\frac{\d g}{\d\bar z^l}-
\frac{\d g}{\d \bar z^j}\frac{\d f}{\d z^l}).
$$

The {\bf quantizing line bundle} $L\rightarrow M=\G\backslash D$ can be defined
as a line bundle such that the curvature of its natural connection is
the K\"ahler form $\omega$ on $M$. Denoting the canonical line bundle by $L$
we see that the potential 1-form corresponding to the natural connection
on $L$ is $\Theta=i\d\ln (s,s)=-i\d\ln K(z,z)$, hence the curvature
$d\Theta=-i\bar \d\d\ln K(z,z)=\omega$ and this is indeed a quantizing line 
bundle for $M$.

A holomorphic function $f:D\rightarrow\C$ is called an {\bf automorphic
form of weight $\bf k$} if 
\begin{equation}
f(\g z)[\det J(\g,z)]^k=f(z)
\label{autlaw}
\end{equation}
for any $z\in D$, $\g\in\G$; here $J(\g,z)$ is the Jacobi matrix
of transformation $\g$ at point $z$. In the context of 1.1 the automorphy
factor $\mu(\g,z)=[\det J(\g,z)]^k$. It is immediately clear that 
automorphic forms of weight $k$ form the complex inner product space
$H^0(M,L^{\otimes k})$ of holomorphic sections of $L^{\otimes k}$.

Now we consider a family of maps $p_k$, here $k$ is a positive integer, 
such that 
$p_k(f)=T_f^{(k)}$, where $f$ belongs to the Poisson algebra of smooth
real-valued functions on $M$ and $T_f^{(k)}$ is the Toeplitz operator
on $H^0(M,L^{\otimes k})$ obtained from multiplication
operator $M_f^{(k)}(g)=fg$ on $L^2(M,L^{\otimes k})$ by the orthogonal
compression to the closed subspace $H^0(M,L^{\otimes k})$, i.e.
$T_f^{(k)}=\Pi^{(k)}\circ M_f^{(k)}\circ \Pi^{(k)}$, where 
$\Pi^{(k)}$ is the orthogonal projection from $L^2(M,L^{\otimes k})$
to $H^0(M,L^{\otimes k})$.    

In the Berezin scheme of quantization for each $\hbar=\frac{1}{k}$ 
we consider the space ${\cal F}_{\hbar}$ of functions holomorphic in $D$
and satisfying (\ref{autlaw}) with the scalar product defined by
$$
(f,g)=const(\hbar)\int_M f(z)\bar g(z)[K(z,z)]^{-\frac{1}{\hbar}}d\mu (z),
$$
where $d\mu (z)=\omega^n$ is the $G$-invariant volume form on $D$ corresponding
to the metric (\ref{metric}). 
It is clear that ${\cal F}_\hbar$ is naturally identified
with  $H^0(M,L^{\otimes k})$. 
For the sake of completeness let us also explain briefly how
the operator $\hat A$ corresponding a classical observable $A=A(z)$, 
is defined. First, we consider an analytic continuation $A(z,w)$ of 
the function $A(z)$ to $D\times D$. The covariant symbol $A(z,z)$ 
of $\hat A$ is defined as the diagonal value of the function
$$
A(z,w)=\frac{\int_M\hat{A}[(K(u,w))^{\frac{1}{\hbar}}] 
(K(z,u))^{\frac{1}{\hbar}}d\mu (u)}
{\int_M (K(u,w))^{\frac{1}{\hbar}}(K(z,u))^{\frac{1}{\hbar}}d\mu (u)},
$$
and 
$$
(\hat Af)(z)=const(\hbar)\int_M A(z,w)f(w)[K(z,w)]^{\frac{1}{\hbar}}
[K(w,w)]^{-\frac{1}{\hbar}}d\mu (w).
$$
So we end up with the algebra $A_\hbar$ of covariant symbols of 
bounded operators acting in ${\cal F}_\hbar$.
The $*$-product in $A_\hbar$ is given by
$$
A_1*A_2 \ (z,z)=const(\hbar)\int_M A_1(z,w)A_2(w,z)
(\frac{K(z,w)K(w,z)}{K(z,z)K(w,w)})^{\frac{1}{\hbar}}d\mu(w).
$$
 
In conclusion let us discuss the Poincar\'e series (\ref{p-s}) and  
(\ref{rp-s}). Consider the unitary representation of $G$
in ${\cal F}_\hbar$ given by the operators  
$$
[\pi^k(g)(f)](z)=[\det J(g^{-1},z)]^{k}f(g^{-1}z).  
$$
It can be regarded as a subrepresentation of the left regular representation
of $G$ in $L^2(\G\backslash G)$.
Fix $f\in {\cal F}_\hbar$, then the set $\{ \pi^k(g)(f)|g\in G\}$ is called
a {\it system of coherent states} (strictly speaking, we should
regard two coherent states $\pi^k(g_1)(f)$ and $\pi^k(g_2)(f)$
as equivalent if $\pi^k(g_1)(f)=e^{i\alpha}\pi^k(g_2)(f)$).
Now it is clear that (\ref{p-s}) and (\ref{rp-s}) are just sums
of coherent states corresponding to $f=q$ and the representation
described above.

\subsection{Comments on the subject of the present paper}

In \cite{FK} and in the present paper we consider holomorphic
automorphic forms on $D=\H^n_\C=SU(n,1)/U(n)$. In \cite{FK} we construct
a set of relative Poincar\'e series generating the graded algebra
of $\C$-valued cusp forms on a finite volume quotient of $D$. 
In the present paper we regard holomorphic $\C$-valued automorphic forms
on $\H^n_\C$ as holomorphic sections
of the line bundle $L^{\otimes k}\rightarrow \Gamma\backslash \H^n_\C$, where 
$L$ is a quantizing line bundle on $\Gamma\backslash \H^n_\C$, 
$k$ is an integer, and $\Gamma$ is a discrete cocompact subgroup of $SU(n,1)$. 
We construct relative Poincar\'e series
associated to loxodromic elements of $\G$ and we address 
an interesting problem which is not resolved for Poincar\'e series in general: 
is it true that these series are not identically zero ? 
We restrict ourselves to the case of complex dimension 2 and 
answer ``yes'' to this question going through the following steps:
1) to each hyperbolic element of $\G$ we associate 
a certain Legendrian submanifold of the unit circle bundle in $L^*$
such that the corresponding  Lagrangian submanifold of 
$\G\backslash SU(2,1)/U(2)$ satisfies a Bohr-Sommerfeld condition, 
2) following the method of \cite{BPU} we  conclude that 
the relative Poincar\'e series associated to hyperbolic elements
are not zero for large values of $k$ (i.e. in semi-classical limit 
$\hbar=\frac{1}{k}\rightarrow 0$).



\section{Preliminaries}

\subsection{Complex hyperbolic space}

Consider the complex hyperbolic space 
$$
\H^n_\C=SU(n,1)/S(U(n)\times U(1))=
\P(\{ z\in \C^{n+1} \ | \ \zz <0\} ) \approx B^n,
$$ 
here $B^n$ is the open unit ball in $\C^n$, $\la  \ . \ , \ . \ \ra$ is the Hermitian
product on $\C^{n+1}$ given by 
$\la z,w\ra =z_1\bar w_1+...+z_n\bar w_n-z_{n+1}\bar w_{n+1}$ for 
$z=\pmatrix{z_1 \cr ... \cr z_{n+1}}\in \C^{n+1}$, 
$w=\pmatrix{w_1 \cr ...  \cr w_{n+1}}\in \C^{n+1}$. 

A non-zero vector $z\in\C^{n+1}$ is called {\it negative (null, positive)}
if the value of $\zz$ is negative (null, positive).

For $z,w\in B^n$ the corresponding vectors in $\C^{n+1}$ are 
$z=\pmatrix{z_1 \cr ... \cr z_n \cr 1}$ and 
$w=\pmatrix{w_1 \cr ... \cr w_n \cr 1}$, and
$\la z,w\ra =z_1\bar w_1+...+z_n\bar w_n-1$.    

The group of biholomorphic automorphisms of $\H^n_\C$ is $PU(n,1)=SU(n,1)/center$.
The group $SU(n,1)$ acts on $B^n$ by fractional-linear transformations: for
$$
\g=\elt
$$
we have:
$$
 \g z=\g (z_1,...,z_n)=(\frac{a_{11}z_1 +\dots + a_{1n}z_n+
b_1}{c_1z_1+\dots + c_nz_n+d},\dots, \frac{a_{n1}z_1 +\dots +a_{nn}z_n+
b_n}{c_1z_1+\dots + c_nz_n+d}).
$$
Notice that $\det J(\g,z)=(c_1z_1+\dots + c_nz_n+d)^{-(n+1)}$,
here $J(\g,z)$ denotes the Jacobi matrix of transformation $\g$ at point $z\in B^n$.

An automorphism is called {\it loxodromic} if it has no fixed points in
$B^n$ and fixes two points in $\partial B^n$. Notice that the fixed points of
the automorphisms correspond to the eigenvectors of the corresponding matrices
in $U(n,1)$. A loxodromic  automorphism is called {\it hyperbolic} if it has a lift
to $U(n,1)$ all of whose eigenvalues are real.

A loxodromic element $\g _0  \in SU(n,1)$ has $n-1$ positive eigenvectors
and two null eigenvectors.

Let $v_1$,...,$v_{n-1}$ be the positive eigenvectors of $\g _0$,
$\tau _1$,...,$\tau _{n-1}$ - the corresponding eigenvalues.
Then $|\tau _j|=1$, $1\le j\le n-1$.

Let $X$, $Y$ be the null eigenvectors of $\g _0$. Then the corresponding
eigenvalues are $\l$ and $\bar{\l}^{-1}$ for some $\l \in \C$, $|\l |>1$.

A loxodromic transformation can always be represented by a matrix in $U(n,1)$ 
with eigenvalues $\tau _1$,...,$\tau _{n-1}$,$\l$,$\l^{-1}$ where $\l\in \R$, $|\l|>1$. 

The geodesic connecting $X$ and $Y$ (it is an arc of a circle
orthogonal to $\partial B^n$ or a diameter) is $\g_0$-invariant
and is called the {\it axis} of $\g_0$. The complex line  
containing $X$ and $Y$ ({\it complex geodesic}) is $\g_0$-invariant too.

\subsection{Automorphic forms and geometry of the quotient}

Let $\G$ be a discrete cocompact subgroup of $SU(n,1)$ such that the quotient
$X:=\G\sm B^n$ is smooth.

A holomorphic function $f:B^n\rightarrow \C$ satisfying the automorphy law
\begin{equation}
f(\g z)(\det J(\g,z))^{k}=f(z)
\label{aut-law1}
\end{equation}
for any $\g\in\G$ is a cusp form of weight $(n+1)k$ for $\G$. 
The corresponding automorphic form on $SU(n,1)$ is given by 
$F(g)=f(g(0))\zeta^k$, where $\zeta=\det J(g,0)$, and 
the automorphy law on the group is 
$$
F(\g g)=F(g)
$$
for any $\g\in\G$. Notice that $\g: \zeta\rightarrow \zeta \det J(\g,z)$ 
for any $\g\in SU(n,1)$. 

We shall denote the space of cusp forms 
of weight $(n+1)k$ for $\G$ on $B^n$ by $S_{(n+1)k}(\G)$ and the corresponding 
space of cusp forms on $SU(n,1)$ by $\tilde S_{(n+1)k}(\G)$. The inner product 
on $S_{(n+1)k}(\G)$ and $\tilde S_{(n+1)k}(\G)$ is given by 
$$
(f,g)=(f(z)\zeta^k,g(z)\zeta^k)=i^n\int_{\G\sm B^n}f\bar g (-\zz)^{(n+1)k}
\frac{dz_1\w d\bar z_1\w ... \w dz_n \w d\bar z_n}{(-\zz)^{n+1}}.
$$

Given a subgroup $\G_0$ of $\G$ and a holomorphic function $q(z)\in L^2(\G\sm B^n)$ 
satisfying (\ref{aut-law1}) for all $\g\in\G_0$, 
the {\it relative Poincar\'e series} for $\G_0$ is defined as
$$
\Theta (z)=\summa q(\g z)(\det J(\g,z))^{k},
$$ 
this series is converges absolutely and uniformly on compact sets and 
belongs to the space $S_{(n+1)k}(\G)$. 

The Bergman kernel for the domain $B^n$ is 
$K(z,w)=\frac{1}{(-\la z,w\ra)^{n+1}}$ and the K\"{a}hler form on $X$ is
$$
\Phi _\kappa =\frac{2\kappa i}{\zz ^2}
(\zz \sum_{j=1}^ndz_j\w d\bar{z}_j-
\la dz,z\ra \w \la z,dz\ra )
$$
where $\kappa $ is a positive constant,
it is an $SU(n,1)$-invariant K\"{a}hler form on $B^n$.
\begin{remark}
We will set $\kappa =\frac{n+1}{2}$.
Then the holomorphic sectional curvature
on $B^n$ is $-\frac{2}{\kappa}=-\frac{4}{n+1}$ and the sectional curvature
is pinched between $-\frac{4}{n+1}$ and $-\frac{1}{n+1}$ (\cite{G}
II.2.2, III.1.5).
\end{remark}
Consider the canonical line bundle $L:=\bigwedge ^nT^*X$ and the dual bundle
$L^*=\bigwedge ^nTX$. The {\it potential 1-form} $\theta$ on $L^*$
is characterized by
$$
\nabla s=-i\theta s,
$$
where $\nabla$ is a connection on $L^*$ and $s$ is the unit section.
The potential 1-form corresponding to the natural connection on $L^*$ is
$$
\theta =i\partial \ln (s,s)=i\partial \ln ((-\zz )^{n+1}).
$$
Here $d=\partial +\bar\partial$.
The {\it curvature} on $L^*$ is $d\theta =-\Phi _{\frac{n+1}{2}}$, 
hence $L$ is the natural quantizing line bundle for $X$.

Let 
$$
E _k=H^0(X,L^{\otimes k})
$$ 
be the complex inner-product
space of holomorphic sections of the $k$-th tensor power of $L$.
Consider the unit circle bundle $P\subset L^*$,
a point of $P$ can be described as $(z,\zeta )$,
where $z\in B^n$ and $\zeta$ is the coordinate on the fiber,
$|\zeta |=(-\zz )^{\frac{n+1}{2}}$. We have:
$$
E _k=\tilde S_{(n+1)k}(\Gamma )=
\{ f(z)\zeta ^k \ | \ (z,\zeta )\in P, \
f \ is \ holomorphic \ on \ B^n, \ 
$$
$$
f(\gamma z)(\det J(\gamma ,z))^{k}=f(z) \
for \ any \ \gamma \in \Gamma \} .
$$

Denote also $\tilde{L}=\bigwedge ^nT^*B^n$,
$\tilde{P}$ - the unit circle bundle in $\tilde{L}^*$.

The {\it connection form} $\alpha :TP\rightarrow \R$ on $P$ is 
$$
\alpha =\theta +i\frac{d\zeta}{\zeta}.
$$
It serves as a contact form on $\tilde{P}$ and $P$.

A Lagrangian submanifold $\L _0\subset X$ satisfies a 
{\it Bohr-Sommerfeld condition} if
$$
\frac{k}{2\pi}\int_C\theta \in \Z
$$
for any closed curve $C\subset \L_0$.
The constant $\frac{1}{k}$ plays role of the Planck constant.

The unit disk bundle $W$ in $L^*$ is a compact, strictly pseudoconvex
domain with smooth boundary. Let us consider the {\it Hardy space} of $P$:
$E \subset L^2(P)$ and the {\it Sz\"ego projector}
$$
\Pi :L^2(P) \rightarrow E
$$
given by the orthogonal projection of $L^2$ onto $E$. We identify:
$$
E =\oplus_{k=0}^\infty E _k.
$$
Let $\tilde E$ be the Hardy space of $\tilde P$,
then $\tilde E=\oplus_{k=0}^\infty \tilde E_k$ where
$$
\tilde E_k=\{ f(z)\zeta ^k \ | \ (z,\zeta )\in \tilde P, \
f \ is \ holomorphic \ on \ B^n, \ 
$$
$$
i^n\int_{B^n} |f(z)|^2(-\zz)^{(n+1)k} 
\frac{dz_1 \w d\bar z_1 \w ... \w dz_n \w d\bar z_n}{(-\zz)^{n+1}}
<\infty \} .
$$
We shall denote the corresponding orthogonal projection by
$$
\tilde\Pi :L^2(\tilde P) \rightarrow \tilde E.
$$

\section{Construction of relative Poincar\'e series associated to certain
loxodromic elements of a discrete cocompact subgroup of SU(n,1)}

Consider a loxodromic automorphism of $B^n$, represented it by a matrix
$\g_0\in U(n,1)$ with eigenvalues 
$\tau_1,...,\tau_{n-1},\l,\l^{-1}$,
$|\tau_j|=1$, $j=1,...,n-1$, $\l\in\R$, $|\l |>1$, denote the corresponding
eigenvectors by $v_1,...,v_{n-1},X,Y$ ($v_1,...,v_{n-1}$ are positive, $X,Y$
are null). Notice that if each $\tau_j$ is a root of 1 then some power
of $\g_0$ is a hyperbolic element.
\begin{assum}
Assume that 1 is among the eigenvalues of $\g_0$. 
\label{eigen1}
\end{assum}
\begin{remark}
If $g\in U(n,1)$ is hyperbolic then $g^2$ is a hyperbolic element of 
$SU(n,1)$ which satisfies Assumption \ref{eigen1} and has the same eigenvectors 
as $g$. 
\end{remark}
Generalizing the construction suggested in \cite{FK},
for any collection, w.l.o.g. $v_1,...v_m$, $m\le n-1$, of positive
eigenvectors corresponding to eigenvalue 1 we construct a relative Poincar\'e series
$$
\Theta_{\g_0,l,k}=\summa q_l(\g z)(\det J(\g,z))^{2k}\in S_{2(n+1)k},
$$
where $\G_0=<\g_0>$,
$$
q_l(z)=\frac{\la z,v_1\ra ^{l_1}...\la z,v_m\ra ^{l_m}}
{(\zx \zy)^{(n+1)k+\frac{l_1+...+l_m}{2}}},
$$
$l_1$,...,$l_m$ are positive integers such that $l_1+...+l_m$ is even,
$l=(l_1,...,l_m)$.
The series converges absolutely in $B^n$ and uniformly on the compact sets
by the Theorem \ref{ct} for $k\ge 1$ (see the Appendix).

In dimension 2 the loxodromic elements of $\G$ satisfying 
Assumption \ref{eigen1} are exactly the hyperbolic elements of $\G$. 
The relative Poincar\'e series associated to a hyperbolic element 
$\g_0\in \G$ is 
$$
\Theta_{\g_0,l,k}=\summa q_l(\g z)(\det J(\g,z))^{2k}\in S_{6k},
$$
where $\G_0=<\g_0>$,
$$
q_l(z)=\frac{\la z,v\ra ^{2l}}
{(\zx \zy)^{3k+l}},
$$
and $l$ is a positive integer.
\begin{remark}
Let $\g_1$ and $\g_2$ be hyperbolic elements of $\G$. 
If $\g_1=\g_2^N$ for a positive integer $N$, then
$\Theta_{\g_1,l,k}=\Theta_{\g_2,l,k}$. If $\g_1$ and $\g_2$
are conjugate in $\G$ then 
$\Theta_{\g_1,l,k}=\Theta_{\g_2,l,k}$.
\end{remark}

\section{Bohr-Sommerfeld tori}

Consider a hyperbolic element $\g _0 \in \G$, denote its null eigenvectors
by $X$, $Y$, denote its positive eigenvector by $v$, then the corresponding
eigenvalues are $\l$, $\l ^{-1}$, $1$, for $\l \in \Bbb R$, $|\l |>1$.

We choose $v$ so that
$$
A=[v \ \ \  \frac{X}{\xy}+\frac{Y}{2} \  \ \ \frac{X}{\xy}-\frac{Y}{2}]
$$
belongs to $SU(2,1)$.

The transformation $A^{-1}$ maps the complex line containing
$X$ and $Y$ to the complex line
$\{ z_1=0 \}$ and maps the geodesic connecting $X$ and $Y$ to the geodesic $\tilde{C}$
connecting $(0,-1)$ and $(0,1)$.
The following loxodromic element of $SU(2,1)$ preserves $\tilde{C}$
and the line $\{ z_1=0  \}$:
$$
\gamma :=A^{-1}\g _0 A=\pmatrix{
1 & 0 & 0 \cr
0 & a & b \cr
0 & b & a
},
$$
$$
a=\frac{\l^2+1}{2\l}, \ b=\frac{\l^2-1}{2\l}, \ a^2-b^2=1.
$$
Denote
$$
w=(w_1,w_2)=A^{-1}z, \ w_1=A^{-1}z_1, \ w_2=A^{-1}z_2,
$$
and apply change of variables
$$
\label{ch-var}
w_2=\frac{re^{i\phi}-i}{re^{i\phi}+i}, \
w_1=\sqrt{1-w_2\bar{w}_2}Re^{i\Theta},
$$
$$
0<\phi <\pi , \ 0<r<+\infty , \ 0<R<1, \ 0\le \Theta <2\pi.
$$
\begin{prop}
Any 2-cylinder $C_{\phi,R}=\{ \phi =const, \ R=const \}$ is $\g$-invariant.
\end{prop}
\begin{remark}
The coordinates $(r,\Theta)$ on $C_{\phi,R}$ are the ``radial'' and the 
``angular'' coordinates respectively.
\end{remark}
\proof
$$
w_2\rightarrow \frac{aw_2+b}{bw_2+a}=
\frac{a\frac{re^{i\phi}-i}{re^{i\phi}+i}+b}
{b\frac{re^{i\phi}-i}{re^{i\phi}+i}+a}=
\frac{r(a+b)e^{i\phi}-i(a-b)}{r(a+b)e^{i\phi}+i(a-b)}=
\frac{r\frac{a+b}{a-b}e^{i\phi}-i}{r\frac{a+b}{a-b}e^{i\phi}+i},
$$
so
$$
r\rightarrow r\frac{a+b}{a-b}, \ \ \phi \rightarrow \phi ,
$$
also
$$
\frac{|w_1|}{\sqrt{1-w_2\bar w_2}}\rightarrow
\frac{|\frac{w_1}{bw_2+a}|}
{\sqrt{1-|\frac{aw_2+b}{bw_2+a}|^2}}=
\frac{|w_1|}{\sqrt{1-w_2\bar w_2}},
$$
so $R\rightarrow R$.

\noindent $\bigcirc$

For a positive integer $l$ consider the following submanifold of $\tilde P$:
$$
\tilde{T}(l):=\{ (w,(-\ww )^{\frac{3}{2}}e^{i\psi}) \ | \
w\in \{ \phi =\frac{\pi}{2}, \ R=\sqrt{\frac{l}{3k+l}} \}, \
\psi =-\frac{l}{k}\Theta \}.
$$
Denote $T(l):=<\g >\setminus \tilde{T}(l)$,
$\tilde \L(l):=A\tilde T(l)$ and $\G_0:=<\g_0>$.
\begin{prop}
$\L (l):=AT(l)=\G_0\setminus \tilde{\L}(l)$
is a compact Legendrian submanifold of $P$.
\end{prop}
\proof
$T(l)$ and $\L (l)$ are compact submanifolds of $P$.

Let us prove that $\L (l)$ is Legendrian.
Submanifolds $T(l)$ and $\L (l)$ have dimension $2$.
The restriction of $\alpha$ onto $\tilde T(l)$ is
$$
-3i\frac{\la dw,w\ra}{\ww}+i\frac{d\zeta}{\zeta}=
-3i\frac{\la dw,w\ra}{\ww}-d\psi +i\frac{d(-\ww )^\frac 32}{(-\ww )^\frac 32}=
$$
$$
-3i\frac{\la dw,w\ra}{\ww}-d\psi -\frac 32
i\frac{\la dw,w\ra +\la w,dw\ra}{(-\ww )^\frac 32}(-\ww )^\frac 12=
$$
$$
-3i\frac{\la dw,w\ra}{\ww}-d\psi
+\frac 32 i\frac{\la dw,w\ra }{\ww}
+\frac 32 i\frac{\la w,dw\ra }{\ww} =
$$
$$
-\frac 32 i\frac{\bar{w}_1dw_1+\bar{w}_2dw_2}{w_1\bar{w}_1+w_2\bar{w}_2-1}
+\frac 32 i\frac{w_1d\bar w_1+w_2d\bar w_2}{w_1\bar{w}_1+w_2\bar{w}_2-1}
-d\psi =
$$
$$
-\frac 32i\frac{(1-w_2\bar{w}_2)R^2e^{-i\Theta}ie^{i\Theta}d\Theta +
\sqrt{1-w_2\bar{w}_2}R^2e^{-i\Theta}e^{i\Theta}d\sqrt{1-w_2\bar{w}_2}
+\bar w_2dw_2}{(1-w_2\bar{w}_2)R^2+w_2\bar{w}_2-1}+
$$
$$
\frac 32i\frac{-(1-w_2\bar{w}_2)R^2e^{-i\Theta}ie^{i\Theta}d\Theta +
\sqrt{1-w_2\bar{w}_2}R^2e^{-i\Theta}e^{i\Theta}d\sqrt{1-w_2\bar{w}_2}
+w_2d\bar w_2}{(1-w_2\bar{w}_2)R^2+w_2\bar{w}_2-1}-d\psi =
$$
$$
-\frac 32i\frac{(1-w_2^2)R^2id\Theta +
\sqrt{1-w_2^2}R^2d\sqrt{1-w_2^2}
+w_2dw_2}{(1-w_2^2)(R^2-1)}+
$$
$$
\frac 32i\frac{-(1-w_2^2)R^2id\Theta +
\sqrt{1-w_2^2}R^2d\sqrt{1-w_2^2}
+w_2dw_2}{(1-w_2^2)(R^2-1)}-d\psi =
$$

$$
-3i\frac{(1-w_2^2)R^2id\Theta}{(1-w_2^2)(R^2-1)}-d\psi =
3\frac{ R^2d\Theta }{ R^2-1 }-d\psi =
3\frac{ \frac{l}{3k+l}d\Theta }{\frac{l}{3k+l}-1}+
\frac{l}{k}d\Theta =0.
$$
We showed that $T(l)$ is a Legendrian submanifold of $P$.
To prove that $AT(l)$ is Legendrian too it is enough to show that
$\alpha$ is $SU(2,1)$-invariant.

Let $M\in SU(2,1)$,
$$
M:(z,\zeta )\rightarrow (Mz,\zeta \det J(M,z))=(Mz,\zeta c^3),
$$
where $c=c(z)=(m_{31}z_1+m_{32}z_2+1)^{-1}$.
$$
\alpha =i\frac{d\zeta}{\zeta} -3i\frac{\la dz,z\ra}{\la z,z\ra}=
i\frac{d\zeta}{\zeta} -3i\partial \ln (-\zz ).
$$
We have:
$$
i\frac{d(c^3\zeta)}{c^3\zeta} -3i\partial \ln (-\la Mz,Mz\ra )=
i\frac{c^3d\zeta +3c^2\zeta dc}{c^3\zeta} -
3i\partial \ln (-\zz c\bar{c})=
$$
$$
i\frac{d\zeta}{\zeta}+3i\frac{dc}{c}-
3i\partial \ln (-\zz )-3i \partial\ln (c\bar{c})=
i\frac{d\zeta}{\zeta}+3i\frac{\partial c}{c}-
3i\partial \ln (-\zz )-3i\frac{\partial c}{c}=
$$
$$
i\frac{d\zeta}{\zeta}-3i\partial \ln (-\zz ).
$$
\noindent $\bigcirc$

The natural projection $\L _0(l)$ of $\L (l)$ onto $X$
is a compact Lagrangian submanifold of $X$.
\begin{prop}
$\L _0(l)$ satisfies a Bohr-Sommerfeld condition.
\end{prop}
\proof
Let $\tilde T_0(l)$ be the natural projection of $\tilde T(l)$ onto $B^2$,
and let $T_0(l)$ be the natural projection of $T(l)$ onto $X$,
$AT_0(l)=\L _0(l)$. If $C\subset \L_0(l)$
is a closed curve then $A^{-1}C\subset T_0(l)$ is also closed.
Let $z\in \L_0(l)$, $w\in T_0(l)$, $c=c(w)=(a_{31}w_1+a_{32}w_2+a_{33})^{-1}$,
we have:
$$
\int_C\theta =-3i\int_C\partial\ln (-\zz )=
-3i\int_{A^{-1}C}\partial\ln (-\la Aw,Aw\ra )=
$$
$$
-3i\int_{A^{-1}C}\partial\ln (-\la w,w\ra c\bar c)=
$$
$$
-3i\int_{A^{-1}C}\partial \ln (-\la w,w\ra ) +\frac{\bar c\partial c}{c\bar c}=
-3i\int_{A^{-1}C}\partial \ln (-\la w,w\ra ) +\partial \ln c=
$$
$$
-3i\int_{A^{-1}C}\partial \ln (-\la w,w\ra ) +d\ln c=
-3i\int_{A^{-1}C}\partial \ln (-\la w,w\ra ),
$$
so $\int_C\theta$ is $A^{-1}$-invariant (in fact $SU(2,1)$-invariant) and it is
enough to prove that $T_0(l)$ satisfies the Bohr-Sommerfeld condition.
The restriction of the potential 1-form onto $\tilde T_0(l)$ is
$$
-3i\frac{\bar w_1dw_1+\bar w_2dw_2}{w_1\bar w_1+w_2\bar w_2-1}=
$$
$$
-3i\frac{(1-w_2\bar{w}_2)R^2e^{-i\Theta}ie^{i\Theta}d\Theta +
\sqrt{1-w_2\bar{w}_2}R^2e^{-i\Theta}e^{i\Theta}d\sqrt{1-w_2\bar{w}_2}
+\bar w_2dw_2}{(1-w_2\bar{w}_2)R^2+w_2\bar{w}_2-1}=
$$
$$
-3i\frac{(1-w_2^2)R^2id\Theta +\sqrt{1-w_2^2}R^2d\sqrt{1-w_2^2}
+w_2dw_2}{(1-w_2^2)(R^2-1)}=
$$
$$
-3i\frac{(1-w_2^2)R^2id\Theta +\sqrt{1-w_2^2}R^2\frac{-2w_2dw_2}{2\sqrt{1-w_2^2}}
+w_2dw_2}{(1-w_2^2)(R^2-1)}=
$$
$$
-3i(\frac{R^2i}{R^2-1}d\Theta -\frac{w_2dw_2}{1-w_2^2})=
$$
$$
3\frac{R^2}{R^2-1}d\Theta -3i\frac 12 d\ln (1-w_2^2)=
-\frac lk d\Theta -\frac 32id\ln (1-w_2^2),
$$
then
$$
\frac{3k}{2\pi}\int_{A^{-1}C}(-\frac lk d\Theta -\frac 32id\ln (1-w_2^2))=
$$
$$
\frac{-3l}{2\pi}\int_{A^{-1}C}d\Theta =
\frac{-3l}{2\pi}2\pi m=-3lm\in \Z .
$$
\noindent $\bigcirc$

So the torus $\L _0(l)$ is a Lagrangian submanifold satisfying the 
Bohr-Sommerfeld condition.

\begin{prop}
The orthogonal projection
of the delta function at $(w,\eta )\in \tilde P$ into $\tilde E_k$
is the function
$$
\Psi _{(w,\eta )}(z,\zeta ):=\tilde\Pi _k(\delta _{(w,\eta )})=
\frac{(3k-1)(3k-2)}{4\pi ^2}
\frac{ \zeta ^k \bar{\eta}^k }{\la z,w\ra ^{3k}}.
$$
\end{prop}
\begin{remark}
The orthogonal projection
of the delta function at $(w,\eta )\in \tilde P$ into $\tilde E_k$
is the {\it coherent state} in $\tilde E_{k}$ associated
to the point $(w,\eta )\in\tilde P$, by definition
$g\Psi _{(w,\eta )}=\Psi _{g(w,\eta )}$ for $g\in SU(2,1)$.
\end{remark}
\proof
The fact that $\Psi _{(w,\eta )}=\tilde\Pi _k(\delta _{(w,\eta )})$ is
equivalent to the reproducing property:
$$
F(w,\eta )=\int_{\tilde P}\bar\Psi_{(w,\eta )}(z,\zeta )F(z,\zeta )dV\w d\psi
$$
for all $F\in \tilde E_k$. Given any orthonormal basis $\{ F_{l,k}\}$ for
$\tilde E_k$,
we can write the reproducing kernel as the series
$\Psi_{(w,\eta )}(z,\zeta )=\sum_l \bar F_{l,k}(w,\eta )F_{l,k}(z,\zeta )$
which converges absolutely and uniformly on compact sets.

Using the basis 
$$
F_{l,m,k}(z,\zeta )=
\frac{1}{2\pi}\sqrt{\frac{(3k+l+m-1)!}{l!m!(3k-3)!}}z_1^lz_2^m\zeta ^k,
$$
which is orthonormal with respect to the inner product
$$
(f(z)\zeta^k,g(z)\zeta^k)=i^2\int_{B^2}f\bar{g}(-\zz )^{3k-3}
dz_1\w d\bar{z}_1\w dz_2\w d\bar{z}_2,
$$
we obtain:
$$
\Psi _{(w,\eta )}(z,\zeta )=
\sum_{l,m} \bar{F}_{l,m,k}(w,\eta )F_{l,m,k}(z,\zeta )=
$$
$$
\sum_{l,m}\frac{1}{(2\pi )^2}\frac{(3k+l+m-1)!}{l!m!(3k-3)!}
\bar{w}_1^l\bar{w}_2^m z_1^lz_2^m \zeta ^k\bar\eta ^k=
$$
$$
\frac{\zeta ^k\bar\eta ^k}{4\pi ^2(3k-3)!}
\sum_{l,m}\frac{(3k+l+m-1)!}{l!m!}(\bar{w}_1z_1)^l(\bar{w}_2z_2)^m.
$$
To calculate
$$
\sum_{l,m}\frac{(3k+l+m-1)!}{l!m!}x^ly^m=
\sum_m \frac{y^m}{m!}\sum_l \frac{(3k+l+m-1)!}{l!}x^l
$$
let us find first $\sum_l\frac{(N+l)!}{l!}t^l$. Integrating once we get
$\sum_{l=0}^{\infty}\frac{(N+l)!}{(l+1)!}t^{l+1}$; integrating twice
we get $\sum_{l=0}^{\infty}\frac{(N+l)!}{(l+2)!}t^{l+2}$; integrating
$N$ times we get:
$$
\sum_{l=0}^{\infty}\frac{(N+l)!}{(l+N)!}t^{l+N}=
\sum_{l=0}^{\infty}t^{l+N}=\frac{t^N}{1-t},
$$
differentiating this expression $N$ times we get:
$$
\sum_l\frac{(N+l)!}{l!}t^l=(\frac{t^N}{1-t})^{(N)}=
(\frac{t^N-1+1}{1-t})^{(N)}=
$$
$$
(\frac{(t-1)(t^{N-1}+t^{N-2}+...+t+1)+1}{1-t})^{(N)}=
(\frac{1}{1-t})^{(N)}=\frac{N!}{(1-t)^{N+1}},
$$
hence
$$
\sum_m \frac{y^m}{m!}\sum_l \frac{(3k+l+m-1)!}{l!}x^l=
\sum_m \frac{y^m}{m!} \frac{(3k+m-1)!}{(1-x)^{3k+m}}=
$$
$$
\frac{1}{(1-x)^{3k}}\sum_m \frac{(3k+m-1)!}{m!}(\frac{y}{1-x})^m=
\frac{1}{(1-x)^{3k}}\frac{(3k-1)!}{(1-\frac{y}{1-x})^{3k}}=
\frac{(3k-1)!}{(1-x-y)^{3k}},
$$
so
$$
\Psi _{(w,\eta )}(z,\zeta )=
\frac{\zeta ^k\bar\eta ^k}{4\pi ^2(3k-3)!}
\frac{(3k-1)!}{(1-\bar{w}_1z_1-\bar{w}_2z_2)^{3k}}=
\frac{(3k-1)(3k-2)}{4\pi ^2}\zeta ^k\bar\eta ^k
\frac{1}{(-\la z,w\ra )^{3k}}.
$$
\noindent $\bigcirc$

For weight $6k$ we have:
$$
\Psi _{(u,\eta )}(z,\zeta )=\tilde\Pi _{2k}(\delta _{(u,\eta )})=
\frac{(6k-1)(6k-2)}{4\pi ^2}
\frac{ \zeta ^{2k} \bar{\eta}^{2k} }{\la z,u\ra ^{6k}}.
$$
We omit the weight in the notation $\Psi _{(u,\eta )}(z,\zeta )$
but further exposition will be for weight $6k$ so this will not lead
to any confusion.

To get the orthogonal projection of the delta function at
$[(u,\eta )]\in P=\G\sm\tilde P$ (by $[(u,\eta )]$ we denote the equivalence
class of $(u,\eta )$) into $E_{2k}$ we average over the action of $\G$:
\begin{equation}
\Pi_{2k}(\delta _{[(u,\eta )]})=\sum_{g\in\G}g\Psi_{(u,\eta )}.
\label{del-proj}
\end{equation}
The function $\Psi_{(u,\eta )}$ belongs to $\tilde E_{2k}$, hence
the series (\ref{del-proj}) converges absolutely and uniformly on compact sets
by Theorem 9.1 \cite{Borel}.

Following the method of \cite{BPU}, to the submanifold
$\Lambda (l)\subset P$ with a half-form $\nu$
we associate a function
$$
\Phi :=\int_{\L (l)}\Pi_{2k}(\delta _{[(u,\eta )]})\nu =
\sum_{g\in\G}\int_{\L (l)}g\Psi_{(u,\eta )}\nu =
$$
$$
\sum_{g\in \G /\G_0}\sum_{m=-\infty}^{+\infty}
\int_{\L (l)}g\g_0^m\Psi_{(u,\eta )}(z,\zeta )\nu =
\sum_{g\in \G /\G_0}\sum_{m=-\infty}^{+\infty}
\int_{\L (l)}g\Psi_{\g_0^m(u,\eta )}(z,\zeta )\nu =
$$
$$
\sum_{g\in \G /\G_0}
g\int_{\tilde\L (l)}\Psi_{(u,\eta )}(z,\zeta )\nu =
\sum_{g^{-1}\in\G_0\sm\G}
g^{-1}\int_{\tilde\L (l)}\Psi_{(u,\eta )}(z,\zeta )\nu .
$$

\begin{prop}
$$
\int _{\tilde\L (l)}\Psi _{(u,\eta )}(z,\zeta )\nu =
C\frac{\la z,v\ra ^{2l} }{(\la z,X\ra \la z,Y\ra )^{3k+l}}\zeta^{2k},
$$
where the constant $C$ is given by
$$
C=2^{3k+l-2}\frac{i}{(2l)!}\frac{(6k+2l-1)!}{(6k-3)!}
\frac{ (3k)^{3k}l^l }{ (3k+l)^{3k+l} }
\la Y,X\ra ^{3k+l}
\sum_{j=0}^{3k+l-1}\frac{(-1)^j(3k+l-1)!}{j!j!(6k+2l-1-j)!}.
$$
and we take
$$
\nu =\frac{d(A^{-1}u_1)}{A^{-1}u_1}\w \frac{d(A^{-1}u_2)}{1-(A^{-1}u_2)^2}.
$$
\end{prop}
\begin{remark}
The half-form $\nu$ on $\tilde \L(l)$ is $\g_0$-invariant and in properly
chosen coordinates $(r,\Theta)$ on $\tilde \L_0(l)$ it is expressed as
$\nu=\frac{i}{2}d\Theta\w\frac{dr}{r}$. 
\end{remark}
\proof
Let $u\in \tilde\L (l)$, $w=A^{-1}u\in \tilde T(l)$, then
$$
\int _{\tilde\L (l)}\Psi _{(u,\eta )}(z,\zeta )\nu =
\frac{(6k-1)(6k-2)}{4\pi ^2}
\int _{\tilde\L (l)}\frac{ \zeta ^{2k} \bar{\eta}^{2k} }{\la z,u\ra ^{6k}}
\frac{d(A^{-1}u_1)}{A^{-1}u_1}\w \frac{d(A^{-1}u_2)}{1-(A^{-1}u_2)^2}=
$$
$$
\frac{(6k-1)(6k-2)}{4\pi ^2}\zeta ^{2k}
\int _{\tilde\L (l)} \frac{\bar{\eta}^{2k}}{\la z,u\ra ^{6k}}
\frac{d(A^{-1}u_1)}{A^{-1}u_1}\w \frac{d(A^{-1}u_2)}{1-(A^{-1}u_2)^2}=
$$
$$
\frac{(6k-1)(6k-2)}{4\pi ^2}\zeta ^{2k}
\int _{\tilde T(l)} \frac{((-\ww )^{\frac{3}{2}} e^{-i\psi} \det \bar J(A,w))^{2k}}
{\la z,Aw\ra ^{6k}} \frac{dw_1}{w_1}\w \frac{dw_2}{1-w_2^2}=
$$
$$
\frac{(6k-1)(6k-2)}{4\pi ^2}\zeta ^{2k}
\int _{\tilde T(l)} \frac{(-\ww )^{3k} e^{-i2k\psi} (\det \bar J(A,w))^{2k}}
{\la A^{-1}z,w\ra ^{6k}(\det \bar J(A,w))^{2k}}
(\det J(A^{-1},z))^{2k} \frac{dw_1}{w_1}\w \frac{dw_2}{1-w_2^2},
$$
let $A^{-1}z=\pmatrix{ v_1 \cr v_2 \cr 1 }$, then we get:
$$
\frac{(6k-1)(6k-2)}{4\pi ^2}\zeta ^{2k}(\det J(A^{-1},z))^{2k}
\int _{\tilde T(l)} \frac{(-\ww )^{3k} e^{-i2k\psi}}
{(v_1\bar{w}_1+v_2\bar{w}_2-1)^{6k}}
\frac{dw_1}{w_1}\w \frac{dw_2}{1-w_2^2},
$$
on $\tilde T(l)$
$$w_2=\frac{r-1}{r+1}, \ w_1=\sqrt{1-w_2^2}Re^{i\Theta}=
\frac{2\sqrt{r}}{r+1}Re^{i\Theta},
$$
$$
-\ww =(1-R^2)(1-w_2^2)=(1-R^2)\frac{4r}{(r+1)^2},
$$
so we have:
$$
\frac{(6k-1)(6k-2)}{4\pi ^2}\zeta ^{2k}(\det J(A^{-1},z))^{2k} (4(1-R^2))^{3k}
$$
$$
\int _{\tilde T(l)} \frac{(\frac{r}{(r+1)^2} )^{3k} e^{-i2k\psi}}
{(v_1\frac{2\sqrt{r}}{r+1}Re^{-i\Theta}+v_2\frac{r-1}{r+1}-1)^{6k}}
\frac i2 d\Theta \w \frac{dr}{r}=
$$
$$
\frac{(6k-1)(6k-2)}{4\pi ^2}\zeta ^{2k}(\det J(A^{-1},z))^{2k}
(4(1-R^2))^{3k} \frac i2
$$
$$
\int _0^\infty dr\int_0^{2\pi}d\Theta \frac{r^{3k-1} e^{i2l\Theta}}
{(v_12\sqrt{r}Re^{-i\Theta}+v_2(r-1)-r-1)^{6k}}.
$$
The integral
$$
\int_{|w|=1}\frac{1}{(Aw+B)^{6k}}\frac{dw}{w^{2l+1}}, \ |\frac BA|>1
$$
is equal to
$$
\frac{2\pi i}{(2l)!}\frac{d^{2l}}{dw^{2l}}\frac{1}{(Aw+B)^{6k}}|_{w=0}=
\frac{2\pi i}{(2l)!}\frac{(6k+2l-1)!}{(6k-1)!}\frac{A^{2l}}{B^{6k+2l}},
$$
Let $w=e^{-i\Theta}$, $A=v_12\sqrt{r}R$, $B=v_2(r-1)-r-1$. Let us check that
$|\frac BA|>1$.
$$
|\frac{v_2(r-1)-r-1}{v_12\sqrt{r}R}|=|\frac{v_2w_2-1}{v_1R\sqrt{1-w_2^2}}|>
\frac{|v_2w_2-1|}{\sqrt{1-v_2\bar v_2}R\sqrt{1-w_2^2}}\ge
$$
$$
\frac{|v_2w_2-1|}{\sqrt{1-v_2\bar v_2}\sqrt{1-w_2^2}}\ge 1
$$
because
$$
0\le |v_2-w_2|^2=(\bar v_2-w_2)(v_2-w_2)=
v_2\bar v_2-\bar v_2w_2-w_2v_2+w_2^2=
$$
$$
-\bar v_2w_2-w_2v_2+v_2\bar v_2w_2^2+1+
v_2\bar v_2+w_2^2-v_2\bar v_2w_2^2-1=
$$
$$
(\bar v_2w_2-1)(v_2w_2-1)-(1-v_2\bar v_2)(1-w_2^2).
$$
We get:
$$
\frac{(6k-1)(6k-2)}{4\pi ^2}\zeta ^{2k}(4(1-R^2))^{3k}
\frac 12 \frac{2\pi i}{(2l)!}\frac{(6k+2l-1)!}{(6k-1)!}
(\det J(A^{-1},z))^{2k}
$$
$$
\int _0^\infty \frac{r^{3k-1+l} v_1^{2l}2^{2l}R^{2l}}
{(v_2(r-1)-r-1)^{6k+2l}}dr=
$$
$$
\frac{\pi i}{(2l)!}
\frac{1}{4\pi ^2}\zeta ^{2k}(4(1-R^2))^{3k}
\frac{(6k+2l-1)!}{(6k-3)!}
\frac{v_1^{2l}2^{2l}R^{2l}}{(v_2-1)^{6k+2l}}
(\det J(A^{-1},z))^{2k}
\int _0^\infty \frac{r^{3k-1+l}}
{(r-\frac{v_2+1}{v_2-1})^{6k+2l}}dr.
$$
The integral $\int_0^{\infty}\frac{r^{3k+l-1}}{(r-a)^{6k+2l}}dr=
-res_{z=a}f(z)$ where
$f(z)=\frac{z^{3k+l-1}\ln z}{(z-a)^{6k+2l}}$ (\cite{SFS} 29.6).
Let us calculate this residue:
$$
res_{z=a}f(z)=\frac{1}{(6k+2l-1)!}
\frac{d^{6k+2l-1}}{dz^{6k+2l-1}}z^{3k+l-1}\ln z|_{z=a}=
$$
$$
\frac{1}{(6k+2l-1)!}
\sum_{j=0}^{6k+2l-1}C_{6k+2l-1}^j (z^{3k+l-1})^{(j)}(\ln z)^{(6k+2l-1-j)}
|_{z=a}=
$$
$$
\frac{1}{(6k+2l-1)!}
\sum_{j=0}^{3k+l-1}C_{6k+2l-1}^j (z^{3k+l-1})^{(j)}(\frac{1}{z})^{(6k+2l-2-j)}
|_{z=a}=
$$
$$
\frac{1}{(6k+2l-1)!}
\sum_{j=0}^{3k+l-1}\frac{(6k+2l-1)!}{j!(6k+2l-1-j)!}
\frac{(3k+l-1)!}{j!} z^{3k+l-1-j}
$$
$$
\frac{(6k+2l-2-j)!(-1)^{6k+2l-2-j}}{z^{6k+2l-1-j}}|_{z=a}=
C_1\frac{1}{a^{3k+l}},
$$
where
$C_1=\sum_{j=0}^{3k+l-1}\frac{(-1)^j(3k+l-1)!}{j!j!(6k+2l-1-j)!}$.
Hence
$$
\int_0^{\infty}\frac{r^{3k+l-1}}{(r-\frac{v_2+1}{v_2-1})^{6k+2l}}dr=
-C_1(\frac{v_2-1}{v_2+1})^{3k+l}
$$
and we get:
$$
-\frac{\pi i}{(2l)!}
\frac{1}{4\pi ^2}\zeta ^{2k}(4(1-R^2))^{3k}
\frac{(6k+2l-1)!}{(6k-3)!}C_1
\frac{v_1^{2l}2^{2l}R^{2l}}{(v_2-1)^{6k+2l}}
(\det J(A^{-1},z))^{2k}
(\frac{v_2-1}{v_2+1})^{3k+l}=
$$
$$
\frac{i}{(2l)!}
\frac{1}{4\pi}\zeta ^{2k}(4(1-R^2))^{3k}
\frac{(6k+2l-1)!}{(6k-3)!}C_1
2^{2l}R^{2l}
(\det J(A^{-1},z))^{2k}
\frac{A^{-1}z_1^{2l}}{((A^{-1}z_2)^2-1)^{3k+l}}=
$$
$$
2^{6k+2l-2}\frac{i}{(2l)!}\frac{(6k+2l-1)!}{(6k-3)!}C_1
(\frac{3k}{3k+l})^{3k}(\frac{l}{3k+l})^{l}
(\frac 12 \la Y,X\ra  )^{-3k-l}
\frac{\la z,v\ra ^{2l} }{(\la z,X\ra \la z,Y\ra )^{3k+l}}\zeta ^{2k}.
$$
\noindent $\bigcirc$

We got:
$$
\Phi (z,\zeta )=C\zeta ^{2k}\sum _{g\in \G _0\setminus \G }
q_l(gz)(\det J(g,z))^{2k}\in E_k,
$$
where
$$
q_l(z)=\frac{\la z,v\ra ^{2l} }{(\la z,X\ra \la z,Y\ra )^{3k+l}}
$$
and the relative Poincar\'{e} series associated to $\L _0(l)$ is
$$
\Theta _{\g _0,l,k}(z):=
C\sum _{g\in \G _0\setminus \G } q_l(gz)(\det J(g,z))^{2k} .
$$

\

From the results of \cite{BPU} (Theorem 3.2, Corollary 3.3) it follows that
$\Theta _{\g _0,l,k}$ is non-vanishing for sufficiently large values of $k$.

\

\appendix

\noindent {\bf APPENDIX.}

We shall prove the following theorem modifying the proof of convergence
of Poincar\'e series contained in \cite{Borel} and \cite{B-book}. 

\begin{Th}
Let $\varphi$ be a function on $G=SU(n,1)$. Assume that

1) $\varphi$ is $Z-finite$ 

2) $\varphi\in L^1(\G_0\sm G)$,

3) $\varphi$ is K-finite on the right 

Let $p_\varphi (x)=\summa \varphi (\g x)$.

Then $p_\varphi$ converges absolutely and uniformly on compact sets.
\label{ct}
\end{Th}
\proof

By Lemma 9.2 \cite{Borel} there exists
$\alpha\in C^\infty _c(G)$ satisfying $\alpha(k^{-1}xk)=\alpha(x)$,
$k\in K$, $x\in G$, such that $\varphi =\varphi *\alpha$. 
such that $U^{-1}=U$, the closure of $U$ is compact,
and $U\supset supp \ \alpha$.
We have:
$$
\varphi (\g x)=(\varphi *\alpha )(\g x)=
\int _G \varphi (\g xy)\alpha(y^{-1})dy=
\int _U \varphi (\g xy)\alpha(y^{-1})dy,
$$
hence
$$
|\varphi (\g x)|\le ||\alpha ||_\infty\int_U|\varphi (\g xy)|dy=
||\alpha ||_\infty \int_{xU} |\varphi (\g y)|dy
$$
Here $||\alpha ||_\infty =\sup _{y\in U} |\alpha (y)|$.

Fix a compact subset $C$ of $G$. We want to prove absolute and uniform
convergence on $C$. The closure of $CU$ is compact. $CU$ is covered
by $N$ copies of a "fundamental domain" of $\G$ in $G$,
$N$ is a positive integer (because $\G$ is discrete).
Denote these domains by $F_1$,...,$F_N$.
By a "fundamental domain" of $\G$ in $G$ I mean a connected set
of representatives of $\G\sm G$.

Let $x\in C$. Then
$$
|\varphi (\g x)|\le ||\alpha ||_\infty \int_{xU} |\varphi (\g y)|dy\le
||\alpha ||_\infty \int_{CU} |\varphi (\g y)|dy
$$
and we get
$$
\summa ||\alpha ||_\infty \int_{CU} |\varphi (\g y)|dy=
||\alpha ||_\infty \summa \int_{CU} |\varphi (\g y)|dy\le
$$
$$
||\alpha ||_\infty \summa
(\int_{F_1} |\varphi (\g y)|dy+...+\int_{F_N} |\varphi (\g y)|dy)=
$$
$$
||\alpha ||_\infty (\summa
\int_{F_1} |\varphi (\g y)|dy+...+\summa \int_{F_N} |\varphi (\g y)|dy)=
$$
$$
N||\alpha ||_\infty \int_{\G_0\sm G}|\varphi (y)|dy<\infty .
$$
So we proved that
$$
|\varphi (\g x)|\le c_\g :=
||\alpha ||_\infty \int_{CU}|\varphi (\g y)|dy
$$
and that the numerical series $\summa c_\g$ converges, hence by Weierstrass
theorem the series $\summa \varphi (\g x)$ converges absolutely and
uniformly on $C$. Q.E.D.


\end{document}